\documentclass[12pt,reqno]{article}
\usepackage[english]{babel}
\usepackage{amsmath,amsthm}
\usepackage{amsfonts}
\usepackage{graphicx}
\usepackage{enumerate}
\usepackage[usenames,dvipsnames]{color}
\usepackage{subfigure}
\usepackage[right,pagewise,displaymath, mathlines]{lineno}
\usepackage{epstopdf}
\usepackage{color}
\usepackage{multirow}

\usepackage{tikz}
\usepackage{schemabloc}

\numberwithin{equation}{section}
\numberwithin{figure}{section}
\numberwithin{table}{section}

\newcommand{\nty}{n \to \infty}

\newcommand{\lr}{\left(}
\newcommand{\lb}{\left[}
\newcommand{\rb}{\right]}
\newcommand{\rr}{\right)}
\newcommand{\lf}{\left\{}
\newcommand{\rf}{\right\}}
\newcommand{\lv}{\left|}
\newcommand{\rv}{\right|}

\newcommand{\kn}{k_n}
\newcommand{\mn}{m_n}
\newcommand{\eel}{\end{lemma}}

\newcommand{\xin}{X_{i:n}}
\newcommand{\cin}{c_{i,n}}
\newcommand{\uin}{U_{i:n}}
\newcommand{\wi}{W_i}
\newcommand{\win}{W_{i:n}}
\newcommand{\xia}{\xi_{\alpha}}

\newcommand{\an}{\alpha_n}
\newcommand{\bn}{\beta_n}
\newcommand{\xib}{\xi_{1-\beta}}

\newcommand{\mau}{M_{\alpha}}

\newcommand{\na}{N_{\alpha}}
\newcommand{\nb}{N_{1-\beta}}
\newcommand{\alp}{\alpha}
\newcommand{\aln}{\alpha_n}
\newcommand{\ben}{\beta_n}
\newcommand{\be}{\beta}
\newcommand{\si}{\sigma}
\newcommand{\tn}{L_{n}}

\newtheorem{theorem}{Theorem}[section]

\newtheorem{lemma}{Lemma}[section]

\newtheorem{corollary}{Corollary}[section]

\begin{document}

\noindent{\small Published version:}\\
\noindent{\emph{Probability and Mathematical  Statistics (Wroclav), }}\\
\noindent{\emph{ Vol. 37, Fasc. 1 (2017), pp.~101-118  (open access)}}\\
\noindent{\emph{ doi: 10.19195/0208-4147.37.1.4 }}\\
\begin{center}

{\sf \LARGE Cram\'{e}r  Type Large  Deviations for Trimmed  L-statistics}\footnote{Research partially
supported by the Russian Foundation for Basic Research (grant RFBR no. SS-2504.2014.1).}

\vspace*{7mm}

{\large Nadezhda Gribkova}\footnote{E-mail: n.gribkova@spbu.ru; nv.gribkova@gmail.com}

\medskip

\textit{Faculty of Mathematics and Mechanics, St.\,Petersburg State University,\\ St.\,Petersburg 199034, Russia}

\end{center}

\bigskip

\begin{quote}
\noindent{\bf Abstract.} {\small
In this paper, we propose a~new approach to the investigation  of asymptotic properties of  trimmed
L-statistics and we apply it to the Cram\'{e}r  type large deviation problem.
Our results can be compared with ones in Callaert et~al.~(1982) -- the first and, as far as we know,~the single 
article, where some results on probabilities of large deviations
for the trimmed $L$-statistics were obtained, but under some strict and unnatural conditions.
Our approach is to approximate the trimmed $L$-statistic by a~non-trimmed $L$-statistic (with smooth weight function) based on Winsorized random variables. Using this  method, we establish the Cram\'{e}r  type large deviation results for the trimmed $L$-statistics under quite mild and natural conditions.}

\medskip


\medskip


\end{quote}

\section{Introduction and main results}
\label{imtro}

 Consider  a~sequence $X_1,X_2,\dots $  of independent identically
distributed  real-valued random variables with distribution function  $F$,  and let
$X_{1:n}\le \dots \le X_{n:n}$ denote the
order statistics corresponding to the first $n$ observations.
Define the trimmed L-statistic by 
\begin{equation}
\label{tn}
L_n=n^{-1}\sum_{i=\kn+1}^{n-\mn}\cin\xin,
\end{equation}
where $c_{i,n}\in \mathbb{R}$, \ $k_n$, $m_n$ are two sequences of integers such that $0\le k_n<n-m_n\le n$. Put
$\alpha_n=k_n/n$, \ $\beta_n=m_n/n$. Throughout this paper, we suppose that
$\alpha_n \to \alpha$, $\beta_n\to \beta$, as  $n\to\infty$, where $0<\alpha<1-\beta<1$, i.e. we focus on the case of heavy trimmed $L$-statistic.

In this paper we investigate Cram\'{e}r  type large deviations, i.e. relative errors in the central limit theorem
for $\tn$. First we note that in the case of non-trimmed $L$-statistic ($k_n=m_n=0$) with the coefficients $\cin$ generated by a smooth weight function the
Cram\'{e}r  type large and moderate deviations  were studied in a~number of papers (see Vandemaele and Veraverbeke~\cite{vv82}, Bentkus and Zitikis~\cite{bz90}, Aleskeviciene~\cite{al91}). In contrast, to the best of our knowledge, there exists a~sole  paper -- Callaert et al.~\cite{cvv82} -- devoted to the large deviations for the trimmed $L$-statistics.
However, the result  in~\cite{cvv82} was obtained under some rigorous and unnatural conditions
imposed  on the underlying distribution $F$ and the weights. The method of proof in Callaert et al.~\cite{cvv82} is based on the following two well-known facts: 

1. The joint distribution
of $\xin$ coincides with the joint distribution of $F^{-1}(G(Z_{i:n}))$, $i=1,\dots,n$, where $G$ is the distribution function of the standard exponential distribution, $Z_{i:n}$ are the order statistics corresponding to a~sample of $n$ independent random variable from the distribution $G$. 
 
2. The order statistics $Z_{i:n}$ are distributed as $\sum_{k=1}^iZ_k/(n-k+1)$, where $Z_k$ -- independent standard  exponential random variables. 

These two facts and the  Taylor expansion together enable one to get an~approximation of $\tn$ by a~sum of weighed  i.i.d. random variables for which some suitable known  result on Cram\'{e}r   type large deviations can be applied. This approach was first implemented  by Bjerve~\cite{bj77} to prove a~Berry-Esseen type result for the $L$-statistics. However, use of this method  requires  excessive smoothness conditions imposed on $F$ and leads to the unnatural and complicated  normalization of the $L$-statistic (cf.~Callaert et al.~\cite{cvv82}).

In this article, we propose another approach to  the investigation of asymptotic properties of the trimmed $L$-statistics
different from that used in Bjerve~\cite{bj77} and  Callaert et al.~\cite{cvv82}.  
Our idea is to approximate the trimmed $L$-statistic by a~non-trimmed $L$-statistic with  weights  generated by a~smooth weight function, where  the approximating $L$-statistic is based on the order statistics corresponding to a~sample of $n$ i.i.d. Winsorized random variables. The asymptotic properties that we are interested in are often
well studied in the case of $L$-statistics with a~smooth weight function and bounded observations, this allows us to obtain a~desired result for the trimmed $L$-statistic by applying a~result of the corresponding type to the approximating non-trimmed $L$-statistic; so  it remains only to  evaluate the remainder in the approximation.
Here,  we  apply  our  method to obtain   a~result on probabilities of large
deviations for the trimmed $L$ -statistics, and we establish it  under mild and natural conditions. This our result  on large deviations can be viewed as a~strengthening of the result from Callaert et al.~\cite{cvv82}.

To conclude this introduction, we adduce  a~brief review of the relevant literature.
The class of $L$-statistics is one of the most commonly used classes in statistical inferences.  We refer to
monographs by David and  Nagaraja~\cite{david_2003}, Serfling~\cite{serfling_90}, Shorack and Wellner~\cite{shorack}, van der Vaart~\cite{vanderv} for an~introduction to the theory and applications of  $L$-statistics.
There is a~vast literature on asymptotic properties  of $L$-statistics.
Since we focus on the case of heavy trimmed $L$-statistics,  we will mention mainly  sources  appropriate to our case. The most significant contribution to the establishment of the central limit theorem  for  (trimmed) $L$-statistics was made by Shorack ~\cite{shor69}-\cite{shor72} and  Stigler~\cite{s69}-\cite{s74}.   Mason and Shorack~\cite{mas_shor_90} obtained the necessary and sufficient conditions for the asymptotic normality of the trimmed $L$-statistics. The Berry -- Esseen type bounds under  different sets of conditions were obtained by Bjerve~\cite{bj77}, Helmers~\cite{helm81}-\cite{helm_e2}, Gribkova \cite{gri}. A~great contribution to the research of  second order asymptotic properties for  $L$-statistic was done  by Helmers~\cite{helm80}-\cite{helm_e2}, who  established the Edgeworth expansions for the (trimmed) L-statistics. 
 In  papers  by Bentkus et al.~\cite{bgz}, Friedfich~\cite{friedrich}, Putter and van Zwet~\cite{pz} and van Zwet~\cite{zwet},  the Berry--Esseen type bounds and  Edgeworth expansions for $L$-statistics were derived as the consequences of the very general results for symmetric statistics established in these papers.
Some interesting results on Chernoff's type  large deviations
(for non-trimmed $L$-statistics with smooth weight function) were obtained by Boistard~\cite{boi}. Recently, Gao and Zhao~\cite{gao_zhao} proposed a~general delta method in the theory of Chernoff's type  large deviations and illustrated it by many examples including M-estimators and L-statistics.
A~survey on the $L$-statistics and some   modern applications of them in the economy and theory of actuarial risks can be found in  Greselin~et~al.~\cite{zit_2009}.

We will now proceed to the statement of our results. Define the left-continuous inverse of $F$:
$F^{-1}(u)= \inf \{ x: F(x) \ge u \}$, \ $0<u\le 1$, \
$F^{-1}(0)=F^{-1}(0^+)$, and let $F_n$, $F_n^{-1}$ denote the
empirical distribution function and its inverse respectively.
Let $J$ be a~function defined in an~open set $I$ such that $[\alpha,1-\beta]\subset I\subseteq(0,1)$.
 We will also consider the trimmed $L$-statistics with coefficients  generated by the weight function $J$
 \begin{equation}
 \label{tn0}
\tn^0=n^{-1}\sum_{i=\kn+1}^{n-\mn}\cin^0\xin=\int_{\aln}^{1-\ben}J(u)F_n^{-1}(u)\,du,
\end{equation}
where $\cin^0=n\int_{(i-1)/n}^{i/n}J(u)\,du$.

To state our results, we will need the following set of assumptions.

\smallskip
\noindent{\bf (i)} {\it $J$ is Lipschitz in $I$, i.e. there exists a~constant $C\ge 0$ such that}
\begin{equation}
\label{LipJ}
|J(u)-J(v)|\le C|u-v|,\quad \forall \ \ u,\,v\in I.
\end{equation}

\smallskip
\noindent{\bf (ii)} {\it  $F^{-1}$  satisfies
a~H\"{o}lder condition of order $0<\varepsilon\le 1$ in some neighborhoods $U_{\alpha}$ and $U_{1-\beta}$
of $\alpha$ and $1-\beta$.}

\smallskip
\noindent{\bf (iii)} {\it $\max(|\alpha_n-\alpha|,\,|\beta_n-\beta|)=O\bigl( n^{-\frac 1{2+\varepsilon}}\bigr)$, where $\varepsilon$ is the
H\"{o}lder index from condition}~{\bf (ii)}.

\smallskip
\noindent{\bf (iv)} {\it with $\varepsilon$ from  conditions} {\bf (ii)}-{\bf (iii)}
$$
\sum_{i=k_n+1}^{n-m_n}|c_{i,n}-\cin^0|=O(n^{\frac 1{2+\varepsilon}}).
$$

Define a~sequence of  centering constants
 \begin{equation}
 \label{mun}
\mu_n=\int_{\aln}^{1-\ben}J(u)F^{-1}(u)\,du.
\end{equation}
Since $\aln\to\alp$,  $\ben\to\be$ as $\nty$, both variables $\tn^0$ and $\mu_n$ are well defined for all sufficiently large $n$.

It is well known  (cf.,~e.g.,~\cite{mas_shor_90},~\cite{s74},~\cite{vanderv}) that when  the inverse $F^{-1}$ is continuous at two points $\alp$ and $1-\be$,  smoothness condition~\eqref{LipJ} implies  the weak convergence to the normal law: $\sqrt{n}(\tn^0-\mu_n)\Rightarrow N(0,\si^2)$, where 
\begin{equation}
\label{sigma}
\si^2=\si^2(J,F)=\int_{\alp}^{1-\be}\int_{\alp}^{1-\be} J(u)J(v)(u\wedge v-uv)\, dF^{-1}(u)\,dF^{-1}(v),
\end{equation}
and  $u\wedge v=\min(u,v)$. Here and in the sequel, we use the convention that $\int_a^b=\int_{[a,b)}$ when integrating with respect to the left
continuous integrator $F^{-1}$. All along the article, we assume $\si>0$.

Define the distribution functions of the normalized $\tn$ and $\tn^0$ respectively
\begin{equation}
\label{dfs}
F_{\tn}(x) =\textbf{P}\{\sqrt{n}(\tn-\mu_n)/\si\le x\},\quad F_{\tn^0}(x) =\textbf{P}\{\sqrt{n}(\tn^0-\mu_n)/\si\le x\}.
\end{equation}

Let $\Phi$ denote the standard normal distribution function. Here is our first result on  Cram\'{e}r  type large deviations for $\tn$.

\begin{theorem}
\label{thm1}
Suppose that $F^{-1}$ satisfies condition {\bf (ii)} for some $0<\varepsilon\le 1$ and  the sequences $\aln$ and $\bn$ satisfy  {\bf (iii)}. In addition, assume that the weights $\cin$  satisfy   {\bf (iv)} for some function $J$ satisfying condition {\bf (i)}.

Then for every  sequence $a_n\to 0$ and each $A>0$
\begin{equation}
\label{thm_1}
\begin{split}
1- F_{\tn}(x) &= [1-\Phi(x)](1+o(1)),\\
F_{\tn}(-x)&=\Phi(-x)(1+o(1)),
\end{split}
\end{equation}
as \,$\nty$, uniformly in the range $-A\le x\le a_n n^{\varepsilon/(2(2+\varepsilon))}$.
\end{theorem}
The proof of our main results is relegated to Section~3. Theorem~\ref{thm1} directly implies the following two corollaries.
\begin{corollary}
\label{cor1}
Suppose that the conditions of Theorem~\ref{thm1} are satisfied with $\varepsilon=1$, i.e. $F^{-1}$ is Lipschitz in
some neighborhoods $U_{\alpha}$ and $U_{1-\beta}$
of $\alpha$ and $1-\beta$.  Then for every sequence $a_n\to 0$ and each $A>0$ relations~\eqref{thm_1} hold true, uniformly in the range $-A\le x\le a_n n^{1/6}$.
\end{corollary}

\begin{corollary}
\label{cor2}
Let $\cin=\cin^0=n\int_{(i-1)/n}^{i/n}J(u)\,du$ \,$(\kn+1 \le i\le n-m_n)$, where $J$ is a~function  satisfying
{\bf (i)}. Furthermore, assume that conditions {\bf (ii)} and {\bf (iii)} hold for some $0<\varepsilon\le 1$.
Then  relations~\eqref{thm_1} with  $\tn=\tn^0$ hold true for every sequence $a_n\to 0$ and each $A>0$, uniformly in the range $-A\le x\le a_n n^{\varepsilon/(2(2+\varepsilon))}$.
\end{corollary}

Theorem~\ref{thm1} can be compared with the~result by~Callaert et al.~\cite{cvv82}, where it was assumed that the derivative  $H'=(F^{-1}\circ G)'$ exists and satisfies a~H\"{o}lder condition of order $0<\varepsilon\le 1$  in some
open set, containing $[G^{-1}(\alpha),G^{-1}(1-\beta)]$, where $G$ is the standard exponential distribution function. Moreover, some unnatural condition was imposed on the weights and $H'$
(cf., conditions~(A2) and~(B), Callaert et al.~\cite{cvv82}). In contrast, we use  the natural scale parameter $\sigma$ -- root of the asymptotic variance
of $\tn$ -- for the normalization, and our smoothness condition {\bf (ii)} for $F^{-1}$ is much weaker than one from  Callaert et al.~\cite{cvv82}.

Our Theorem~\ref{thm1} is also related with previous results by Vandemaele and Veraverbeke~\cite{vv82} and  Bentkus and  Zitikis~\cite{bz90} on Cram\'{e}r  type large deviations for non-trimmed $L$-statistics with smooth weight function. The method of proof in the first of these articles was based on Helmers's~\cite{helm81}-\cite{helm_e2} $U$-statistic approximation, and in the second one the $\omega^2$-von Mises statistic type  approximation was applied.
We approximate our trimmed $L$-statistic by $L$-statistics with smooth weight function.
Moreover, we apply the results from the papers mentioned to our approximating non-trimmed $L$-statistic when proving Theorem~\ref{thm1}. Note also that  Cram\'{e}r 's moment conditions
for the underlying distribution assumed in the cited papers  are not needed in the case of the trimmed  $L$-statistics, whereas  the smoothness of $F^{-1}$ near $\alp$ and $1-\be$ becomes essential for the Cram\'{e}r  type large deviations results.

Finally, we state  a~version of Theorem~\ref{thm1}, where the scale factor $\sigma/n^{1/2}$ is replaced by $\sqrt{{\text{\em Var}}(\tn)}$, it is parallel to Theorem~2~(ii) by Vandemaele and Veraverbeke~\cite{vv82}, but now for the trimmed $L$-statistics.

We will need the following two somewhat stronger versions of conditions {\bf (iii)} and {\bf (iv)}.

\smallskip
\noindent{\bf (iii')} \ {\it $\max(|\alpha_n-\alpha|,\,|\beta_n-\beta|)=O\Bigl( n^{-\frac{1}{2+\varepsilon}\lb 1 + \frac{\varepsilon(1-\varepsilon)}{2}\rb}(\log n)^{-\frac{\varepsilon}{2}}\Bigr)$, where $\varepsilon$ is the
H\"{o}lder index from condition}~{\bf (ii)}.

\smallskip
\noindent{\bf (iv')} {\it with $\varepsilon$ from  conditions} {\bf (ii)}-{\bf (iii')}
$$
\sum_{i=k_n+1}^{n-m_n}|c_{i,n}-\cin^0|=O\Bigl( n^{\frac 1{2+\varepsilon}\lb1-\frac{\varepsilon}{2}\rb}\Bigr).
$$
\begin{theorem}
\label{thm2}
Suppose that the conditions of Theorem~\ref{thm1} are satisfied, where {\bf (iii)} and {\bf (iv)} are replaced by   {\bf (iii')} and {\bf (iv')} respectively. \ \ In addition, assume  that \
${\text{\em Var}}(\tn)<\infty$ for all sufficiently large $n$. Then
\begin{equation}
\label{thm_2}
n \sigma^{-2}{\text{\em Var}}(\tn)=1 + O\bigl(n^{-\frac{\varepsilon}{2+\varepsilon}}\bigr).
\end{equation}
Furthermore,  relations~\eqref{thm_1}, where  $\sigma/n^{1/2}$ is replaced by $\sqrt{{\text{\em Var}}(\tn)}$, hold true for every  sequence $a_n\to 0$ and each $A>0$
as \,$\nty$, uniformly in the range $-A\le x\le a_n n^{\varepsilon/(2(2+\varepsilon))}$.
\end{theorem}
Note that  in the  case of heavy trimmed $L$-statistics the condition $\textbf{E}|X_1|^{\gamma}<\infty$ (for some $\gamma>0$) is sufficient for  the finiteness of ${\text{\em Var}}(\tn)$ when  $n$ gets large.
\section{Our method (representation for $\tn^0$ by a~non-trimmed L-statistic)}
\label{lemmas}
Let $\xi_{\nu}=F^{-1}(\nu)$, $0<\nu<1$, be the $\nu$-th quantile of $F$ and  $\wi$ denote $X_i$ Winsorized outside of $(\xia,\xib]$. In other words
\begin{equation}
\label{2_2}
\wi=\left\{
\begin{array}{ll}
\xia,& X_i\le \xia, \\
X_i,& \xia < X_i \le \xib,\\
\xib,& \xib < X_i .
\end{array}
\right.
\end{equation}
Let $\win$ denote the order statistics, corresponding to $W_1,\dots,W_n$ (the sample of $n$ i.i.d. auxiliary random variables).

Define the distribution function $G(x)=\textbf{P}\{\wi\le x\}$ of $\wi$, the corresponding quantile function is equal to
$G^{-1}(u)= \xia \vee (F^{-1}(u) \wedge \xib)$. Here and further on $(a\vee b)=\max(a,b)$. Let $G_n$ and $G_n^{-1}$ denote the corresponding empirical distribution function and its inversion respectively.

We will approximate  $\tn$ by a~linear combination of the order statistics $\win$ with coefficients, generated by the  weight function
\begin{equation}
\label{2_3}
J_w(u)=\left\{
\begin{array}{ll}
J(\alp),& u\le \alp, \\
J(u),& \alp < u \le 1-\be,\\
J(1-\be),& 1-\be < u,
\end{array}
\right.
\end{equation}
which is defined in $[0,1]$. It is obvious that when $J$ is Lipschitz in $I$, i.e. satisfies condition \eqref{LipJ} with some positive constant~$C$, the function $J_w$ is Lipschitz in $[0,1]$ with some constant $C_w\le C$.

Consider the auxiliary non-truncated $L$-statistic given by
 \begin{equation}
 \label{Ln}
\widetilde{L}_n=n^{-1}\sum_{i=1}^{n}\widetilde{c}_{i,n}\win=\int_0^1 J_w(u)G_n^{-1}(u)\,du,
\end{equation}
where $\widetilde{c}_{i,n}=n\int_{(i-1)/n}^{i/n}J_w(u)\,du$. Define the centering constants
 \begin{equation}
 \label{muL}
\mu_{\widetilde{L}_n}=\int_0^1 J_w(u)G^{-1}(u)\,du.
\end{equation}

Since $\wi$ has the finite moments of any order and  because $J_w$ is Lipschitz, the distribution of the normalized  $\widetilde{L}_n$ tends to the standard normal law (see, e.g.,~\cite{s74})
\begin{equation*}
 \label{Ln_to}
\sqrt{n}(\widetilde{L}_n-\mu_{\widetilde{L}_n})/\sigma(J_w,G) \Rightarrow N(0,1),
\end{equation*}
where the asymptotic variance  
\begin{equation}
\label{sigma}
\sigma^2(J_w,G)=\int_0^1\int_0^1 J_w(u)J_w(v)(u\wedge v-uv)\, dG^{-1}(u)\,dG^{-1}(v).
\end{equation}

Observe that for $u\in(\alp,1-\be]$ we have $J_w(u)=J(u)$, $G^{-1}(u)=F^{-1}(u)$, and that $dG^{-1}(u)\equiv 0$ for $u\notin(\alp,1-\be]$. This yields the equality of the asymptotic variances
\begin{equation}
\label{sigma_eq}
\sigma^2(J_w,G)=\sigma^2(J,F)=\sigma^2
\end{equation}
of the truncated $L$-statistic $\tn^0$ and the non-truncated $L$-statistic $\widetilde{L}_n$ based on the Winsorized random variables.

Define the binomial random variable $N_{\nu}= \sharp \{i : X_i \le
\xi_{\nu} \}$, where $0<\nu <1$. Our representation for $\tn^0$ is based on the following simple observation: we see that
\begin{equation}
\label{observ}
\win=\left\{
\begin{array}{ll}
\xia,& i\le \na, \\
\xin,& \na < i \le \nb,\\
\xib,& i> \nb.
\end{array}
\right.
\end{equation}

Put $A_n=\na/n$, \ $B_n=(n-\nb)/n$.
 The following lemma provides us a~useful representation which is crucial in the proof of our main results.
\begin{lemma}
\label{lem_2.1}
\begin{equation}
\label{lem_2.1_1}
\tn^0-\mu_n=\widetilde{L}_n-\mu_{\widetilde{L}_n}+R_n,
\end{equation}
where $R_n=R_n^{(1)}+R_n^{(2)}$,
\begin{equation}
\label{lem_2.1_2}
R_n^{(1)}=\int_{\alp}^{A_n} J_w(u)[F_n^{-1}(u)-\xia]\,du-\int_{1-\be}^{1-B_n} J_w(u)[F_n^{-1}(u)-\xib]\,du
\end{equation}
and
\begin{equation}
\label{lem_2.1_3}
 \ \ \ R_n^{(2)}=\int_{\an}^{\alp} J(u)[F_n^{-1}(u)-F^{-1}(u)]\,du-\int_{1-\bn}^{1-\be} J(u)[F_n^{-1}(u)-F^{-1}(u)]\,du.
\end{equation}
\end{lemma}
\noindent{\bf Proof.} \
First, consider the difference between the centering constants. We obtain
\begin{equation}
\label{ptf_lem_1}
\begin{split}
\mu_{\widetilde{L}_n}-\mu_n=&\int_0^1 J_w(u)G^{-1}(u)\,du - \int_{\aln}^{1-\ben}J(u)F^{-1}(u)\,du=\alp J(\alp)\xia \\
&+\be J(1-\be)\xib -\int_{\an}^{\alp} J(u)F^{-1}(u)\,du + \int_{1-\bn}^{1-\be} J(u)F^{-1}(u)\,du .
\end{split}
\end{equation}
For the difference between $\tn^0$ and $\widetilde{L}_n$ after some simple computations we get
\begin{equation}
\label{ptf_lem_2}
\begin{split}
\tn^0-\widetilde{L}_n=&\int_{\alp}^{1-\be}J(u)[F_n^{-1}(u)-G_n^{-1}(u)]\,du\\
 &+ \int_{\an}^{\alp} J(u)F_n^{-1}(u)\,du- \int_{1-\bn}^{1-\be} J(u)F_n^{-1}(u)\,du\\ &-J(\alp)\int_0^{\alp}G_n^{-1}(u)\,du
 -J(1-\be)\int_{1-\be}^1 G_n^{-1}(u)\,du.
\end{split}
\end{equation}

 Relations~\eqref{ptf_lem_1} and \eqref{ptf_lem_2} together imply
 \begin{equation}
\label{ptf_lem_3}
\tn^0-\widetilde{L}_n +(\mu_{\widetilde{L}_n}-\mu_n)= D_n+ R_n^{(2)},
\end{equation}
where
\begin{equation*}
\label{ptf_lem_4}
\begin{split}
D_n := &\int_{\alp}^{1-\be}J(u)[F_n^{-1}(u)-G_n^{-1}(u)]\,du+ \\
&J(\alp)\lb  \alp\,\xia - \int_0^{\alp}G_n^{-1}(u)\,du\rb + J(1-\be) \lb \be\,\xib- \int_{1-\be}^1 G_n^{-1}(u)\,du\rb.
\end{split}
\end{equation*}

It remains to show that $D_n=R_n^{(1)}$. Let us consider three of six possible cases (treatment for the three other cases is similar and therefore omitted). We use the fact that $F_n^{-1}(u)=G_n^{-1}(u)$ for $A_n<u\le 1-B_n$, \  $G_n^{-1}(u)=\xia$ for $u \le A_n$ and $G_n^{-1}(u)=\xib$ \ for $u>1-B_n$.\\
\noindent{\bf Case 1}.  $\alp\le A_n \le 1-B_n<1-\be$. In this case the second and third terms of $D_n$ are equal to zero, and the first one yields
\begin{equation}
\label{ptf_lem_5}
\begin{split}
D_n =\int_{\alp}^{A_n}J(u)[F_n^{-1}(u)-\xia]\,du+   \int_{1-B_n}^{1-\be} J(u) [F_n^{-1}(u)-\xib]\,du,
\end{split}
\end{equation}
and since $J(u)=J_w(u)$ for $\alp<u\le 1-\be$,  we obtain the desired equality.\\
\noindent{\bf Case 2}.  $\alp\le A_n\le 1-\be < 1-B_n$. In this case we have
\begin{equation}
\label{ptf_lem_6}
\begin{split}
D_n =& \int_{\alp}^{A_n}J(u)[F_n^{-1}(u)-\xia]\,du\\
+ &J(1-\be) \lb \be\,\xib- \int_{1-\be}^{1-B_n} F_n^{-1}(u)\,du -B_n\,\xib\rb\\
=&\int_{\alp}^{A_n}J(u)[F_n^{-1}(u)-\xia]\,du-  \int_{1-\be}^{1-B_n}J(1-\be) [F_n^{-1}(u)-\xib]\,du,
\end{split}
\end{equation}
and since $J(u)=J_w(u)$ for $\alp<u\le A_n$ and $J(1-\be)=J_w(u)$ for $u>1-\be$, the expression on the r.h.s. in \eqref{ptf_lem_5} is equal to $R_n^{(1)}$. \\
\noindent{\bf Case 3}.  $1-\be\le A_n\le 1-B_n$. In this case $D_n$ can be written as
\begin{equation}
\label{ptf_lem_6}
\begin{split}
& \int_{\alp}^{A_n}J_w(u)[F_n^{-1}(u)-\xia]\,du\\
 &- J(1-\be)\int_{1-\be}^{A_n} F_n^{-1}(u)\,du + J(1-\be)\xia(A_n-(1-\be))\\
&+J(1-\be) \lb \be\,\xib-\xia(A_n-(1-\be))- \int_{A_n}^{1-B_n} F_n^{-1}(u)\,du -B_n\,\xib\rb\\
=&\int_{\alp}^{A_n}J_w(u)[F_n^{-1}(u)-\xia]\,du-  \int_{1-\be}^{1-B_n}J(1-\be) [F_n^{-1}(u)-\xib]\,du =R_n^{(1)}.
\end{split}
\end{equation}

This completes the proof of representation~\eqref{lem_2.1_1}. The lemma is proved. \qed

\medskip
In conclusion of this section, we note that the idea of the $L$-statistic approximation emerged as a~result of the   observation of the fact that the asymptotic variances of $\tn^0$ and of the non-trimmed $L$-statistic $\widetilde{L}_n$ based on the Winsorized random variables coincide.
This idea of  $L$-statistic approximation can also be regarded as an~extension of the one used in Gribkova and Helmers~\cite{gh2006}-\cite{gh2007} and~\cite{gh2014} (where the second order asymptotic properties -- the Berry--Esseen bounds and  Edgeworth type expansions -- were established for (slightly) trimmed means and their studentized versions) to the case of trimmed $L$-statistics. In  the papers mentioned,  we constructed  the $U$-statistic type approximations for (slightly) trimmed means using sums of i.i.d. Winsorized observations as the linear $U$-statistic terms; in order to get the quadratic terms, we applied  some special Bahadur--Kiefer representations of von Mises  statistic type for (intermediate) sample quantiles (cf.~Gribkova and Helmers~\cite{gh2011}).

\section{Proof of Theorems~\ref{thm1} and~\ref{thm2}  }
\label{proof}
\noindent{\bf Proof of Theorem~\ref{thm1}}.
Obviously, it suffices to prove the first of relations~\eqref{thm_1}. Set
\begin{equation}
\label{2_1}
V_n=\tn-\tn^0=n^{-1}\sum_{i=\kn+1}^{n-\mn}(\cin -\cin^0)\xin.
\end{equation}
Lemma~\ref{lem_2.1} and relation~\eqref{2_1} together yield
\begin{equation}
\label{proof_1}
\tn-\mu_n=\widetilde{L}_n-\mu_{\widetilde{L}_n}+R_n +V_n.
\end{equation}
In view of the classical Slutsky argument applied to~\eqref{proof_1}, $1- F_{\tn}(x)$ is  bounded above and below by
\begin{equation}
\label{proof_2}
\textbf{P}\{\sqrt{n}(\widetilde{L}_n-\mu_{\widetilde{L}_n})/\si> x-2\delta \}+\textbf{P}\{\sqrt{n}|R_n|/\si> \delta\}+\textbf{P}\{\sqrt{n}|V_n|/\si> \delta\}
\end{equation}
and
\begin{equation}
\label{proof_3}
\textbf{P}\{\sqrt{n}(\widetilde{L}_n-\mu_{\widetilde{L}_n})/\si> x+2\delta \}-\textbf{P}\{\sqrt{n}|R_n|/\si> \delta\}-\textbf{P}\{\sqrt{n}|V_n|/\si> \delta\}
\end{equation}
respectively, for each  $\delta>0$.

Let $z_n=n^{\varepsilon/(2(2+\varepsilon))}$. Fix an~arbitrary sequence $a_n\to 0$ and $A>0$. Without loss of generality we may assume that $a_n\ge 1/\log(1+n)$ (otherwise, we may replace $a_n$ by the~new sequence $a_n'=\max(a_n ,\, 1/\log(1+n))\ge a_n$ without affecting result). Set $\delta=\delta_n=a_n^{-1/2}/z_n$.
From~\eqref{proof_2} and~\eqref{proof_3} it immediately follows  that to prove our theorem it suffices to show that
\begin{equation}
\label{proof_4}
\textbf{P}\{\sqrt{n}(\widetilde{L}_n-\mu_{\widetilde{L}_n})/\si> x \pm 2\delta\}=[1-\Phi(x)](1+o(1)) ,
\end{equation}
\begin{equation}
\label{proof_5}
 \quad \quad \ \ \,\textbf{P}\{\sqrt{n}|R_n|/\si> \delta\}=[1-\Phi(x)]o(1) ,
\end{equation}
\begin{equation}
\label{proof_6}
\quad \quad \ \ \, \textbf{P}\{\sqrt{n}|V_n|/\si> \delta\}=[1-\Phi(x)]o(1) ,
\end{equation}
uniformly in the range $-A\le x\le a_n z_n$.

\noindent{\bf Proof of \eqref{proof_4}}. Since $\widetilde{L}_n$ is the non-truncated linear combination of order statistics corresponding to the sample $W_1,\dots,W_n$ of i.i.d. bounded random variables and because its weight function $J_w$ is  Lipschitz in $[0,1]$, we can apply the results on probabilities of large deviations by  Vandemaele and Veraverbeke~\cite{vv82} and  by Bentkus and Zitikis~\cite{bz90}. Set $B=A+2\sup_{n\ge 1}\delta_n$ and  $b_n=a_n+2\delta_n$. Since $a_n\ge 1/\log(1+n)$, the number $B$ exists, and $b_n\to 0$. Then, by   Theorem~2~({\em i}), of Vandemaele and Veraverbeke~\cite{vv82}  for  $x:$  $-B\le x \pm 2\delta <0$, and by Theorem~1.1 of Bentkus and Zitikis~\cite{bz90} for $x:$  $0 \le x \pm 2\delta \le b_n n^{1/6}$), we obtain
\begin{equation}
\label{proof_7}
\textbf{P}\{\sqrt{n}(\widetilde{L}_n-\mu_{\widetilde{L}_n})/\si> x \pm 2\delta \}=[1-\Phi(x \pm 2\delta)](1+o(1)) ,
\end{equation}
uniformly with respect to $x$ such that $-B\le x \pm 2\delta \le b_n n^{1/6}$. In particular, relation~\eqref{proof_7} holds true
uniformly in the range $-A\le x\le a_n n^{1/6}$.  To prove~\eqref{proof_4}, it remains to note that
since $2\,\delta a_n z_n =2\sqrt{a_n}\to 0$, Lemma~A.1 from Vandemaele and Veraverbeke~\cite{vv82}  now yields
\begin{equation}
\label{proof_8}
1-\Phi(x \pm 2\delta)=[1-\Phi(x)](1+o(1)),
\end{equation}
as $\nty$, uniformly in the  range $-A\le x\le a_n z_n$.

\noindent{\bf Proof of \eqref{proof_5}}. Let $I_1^{(j)}$ and $I_2^{(j)}$ denote the first and the second terms of $R^{(j)}_n$ (cf.~\eqref{lem_2.1_2}--\eqref{lem_2.1_3}) respectively, $j=1,2$. In this notation,    $R_n=I_1^{(1)}-I_2^{(1)}+I_1^{(2)}-I_2^{(2)}$ and
\begin{equation}
\label{proof_9}
\textbf{P}\{\sqrt{n}|R_n|/\si> \delta\}\le \sum_{k=1}^2 \textbf{P}\{\sqrt{n}|I^{(1)}_k|/\si> \delta/4\}
+\sum_{k=1}^2 \textbf{P}\{\sqrt{n}|I^{(2)}_k|/\si> \delta/4\} .
\end{equation}
Thus, it suffices to show that for each positive  $C$ (in particular, for $C=\si/4$),
\begin{equation}
\label{proof_10}
\textbf{P}\{\sqrt{n}|I^{(j)}_k|>C \delta \}=[1-\Phi(x)]o(1) , \text{\ \ $k,j=1,2,$}
\end{equation}
as $\nty$, uniformly in the range $-A\le x\le a_n z_n$. We will prove~\eqref{proof_10} for $I^{(1)}_1$ and $I^{(2)}_1$ (the treatment of $I^{(1)}_2$ and $I^{(2)}_2$ is similar and therefore omitted).

Consider $I^{(1)}_1$. First, note that if $\alp < A_n$, then $\max_{u\in (\alp,A_n)} |F_n^{-1}(u)-\xia|=\xia -X_{[n\alp]+1:n}\le \xia -X_{[n\alp]:n}$, as $F_n^{-1}$ is monotonic. Here and in what follows $[x]$ represents the greatest integer function. Similarly we find that if $A_n \le \alp$, then $\max_{u\in (A_n,\alp)} |F_n^{-1}(u)-\xia|=X_{[n\alp]:n}-\xia$. Furthermore, by the Lipschitz condition for $J$, there exists a~positive $K$ such that $\max_{u\in [0,1]}J_w(u)\le \sup_{u\in I}J(u)\le K$. This yields
\begin{equation}
\label{proof_11}
|I^{(1)}_1| =\lv \int_{\alp}^{A_n} J_w(u)[ F_n^{-1}(u)-\xia]\, du \rv  \le K |A_n-\alp| |X_{[n\alp]:n}-\xia|.
\end{equation}
Define a~sequence of intervals $\Gamma_n=[\alp\wedge\an,\alp\vee\an+1/n)$,  then we obtain
\begin{equation}
\label{proof_11_}
|I^{(2)}_1| =\lv \int_{\an}^{\alp}J(u) [F_n^{-1}(u)-F^{-1}(u)]\, du \rv \le K |\an-\alp| \emph{D}_n,
\end{equation}
where $\emph{D}_n= \max_{i:\,i/n\in \Gamma_n}|\xin-F^{-1}(i/n)|\vee |\xin-F^{-1}((i-1)/n)|$.

Let $U_1,\dots,U_n$ be a~sample of independent $(0,1)$-uniform distributed random variables, $\uin$ -- the corresponding order statistics. Set $\mau=\sharp \{i : U_i \le \alp \}$.  Since the joint
distribution of  $\xin$ and $\na$ coincides with the joint distribution of $F^{-1}(\uin)$ and $\mau$, $i=1,\dots,n$, in order to prove~\eqref{proof_10}, it suffices to show  that
\begin{equation}
\label{proof_12}
\begin{split}
\textbf{P}\{|\mau-n\alp| |U_{[n \alp ]:n}-\alp|^{\varepsilon} > C\sqrt{n}\,\delta \}&=[1-\Phi(x)]o(1) ,\\
\textbf{P}\{\sqrt{n}|\an-\alp| \emph{D}_{n,u}^{\varepsilon} > C\,\delta \}&=[1-\Phi(x)]o(1) ,\\
\textbf{P}\Bigl( \bigcup_{i:\,i/n\in \Gamma_n}\lf \uin \notin U_{\alp}\rf) &=[1-\Phi(x)]o(1) ,
\end{split}
\end{equation}
as $\nty$, uniformly in the range $-A\le x\le a_n z_n$. Here $U_{\alp}$ is the~neighborhood of $\alp$ , in which  $F^{-1}$ satisfies a~H\"{o}lder condition of order $\varepsilon$ (cf.~condition~{\bf (ii)}),
\begin{equation}
\label{proof_12a}
\emph{D}_{n,u}^{\varepsilon}=\max_{i:\,i/n\in \Gamma_n}|\uin-i/n|^{\varepsilon}\vee |\uin-(i-1)/n|^{\varepsilon},
\end{equation}
and $C$ stands for  a~positive constant independent of $n$, which may change  its value,  from line to line.

To shorten notation, let $k=[n\alp]$. Consider the probability on the l.h.s. in the first line of~\eqref{proof_12}. It is equal to
\begin{equation}
\label{proof_13}
 \textbf{P}\{|\mau-n\alp| |U_{k:n}-\alp|^{\varepsilon}> C  a_n^{-\frac 12} n^{\frac 1{2+\varepsilon}} \}
\le  \textbf{P}_1+\textbf{P}_2,
\end{equation}
where
\begin{equation*}
\begin{split}
&\textbf{P}_1:=\textbf{P}\{|\mau-n\alp| > C_1 a_n^{-\frac 1{2(1+\varepsilon)}} n^{\frac{1+\varepsilon}{2+\varepsilon}} \}, \\
&\textbf{P}_2:=\textbf{P}\{|U_{k:n}-\alp|^{\varepsilon}> C_2 a_n^{-\frac {\varepsilon}{2(1+\varepsilon)}} n^{-\frac{\varepsilon}{2+\varepsilon}} \},
\end{split}
\end{equation*}
$C_1$, $C_2$ are any positive constants such that $C_1C_2=C$. Let us estimate $\textbf{P}_1$ and  $\textbf{P}_2$. Set  $h=C_1a_n^{-\frac 1{2(1+\varepsilon)}} n^{\frac{1+\varepsilon}{2+\varepsilon}-1}$. Since $h<1-\alp$ for all sufficiently large $n$ (because $a_n\ge 1/\log(1+n)$), by Theorem~1 of  Hoeffding~\cite{ho} we have
\begin{equation}
\label{proof_14}
\textbf{P}_1 = \textbf{P}\{|\mau-n\alp| > n C_1 h\}\le 2\exp(-2nh^2)=2\exp(-2C_1^2n^{\frac {\varepsilon}{2+\varepsilon}}a_n^{-\frac 1{1+\varepsilon}}).
\end{equation}

Next, we evaluate  $1/(1-\Phi(x))$. Let $\phi=\Phi'$. Since  $1-\Phi(x) \sim \phi(x)/x$   as $x\to\infty$, for  $x$ such that $-A\le x\le a_n z_n$ we have
\begin{equation}
\label{proof_phi}
\frac 1{1-\Phi(x)}\le \frac 1{1-\Phi(a_n z_n)} \sim \frac {a_n z_n}{\phi(a_n z_n)}  =  \sqrt{2\pi}\, a_n \, n^{\frac{\varepsilon}{2(2+\varepsilon)}} \exp\lr a_n^2 n^{\frac{\varepsilon}{2+\varepsilon}}/2\rr,
\end{equation}
and combining \eqref{proof_14} and \eqref{proof_phi}, we obtain that
\begin{equation}
\label{proof_15}
\textbf{P}_1 = [1-\Phi(x)]o(1), \ \text{ as} \ n\to\infty,
\end{equation}
uniformly in the range $-A\le x\le a_n z_n$.

Set $p_k=k/({n+1})$, and note that $0< \alp-p_k<n^{-1}$. Then for $\textbf{P}_2$ we have
\begin{equation}
\label{proof_16}
\begin{split}
\textbf{P}_2 \le &\textbf{P}\{|U_{k:n}-p_k|> C_2^{1/\varepsilon} a_n^{-\frac 1{2(1+\varepsilon)}} n^{-\frac 1{2+\varepsilon}} - n^{-1}\},\\
=&\textbf{P}\{\sqrt{n}|U_{k:n}-p_k|> C_2^{1/\varepsilon} a_n^{-\frac 1{2(1+\varepsilon)}} n^{\frac {\varepsilon}{2(2+\varepsilon)}}-n^{-1/2} \}.\\
\end{split}
\end{equation}
Note that the term $n^{-1/2}$ on the r.h.s. in~\eqref{proof_16} is of negligible order and therefore we may omit it.
Set  $\lambda: =C_2^{1/\varepsilon} a_n^{-\frac 1{2(1+\varepsilon)}} n^{\frac {\varepsilon}{2(2+\varepsilon)}} $.
We observe that $\lambda /\sqrt{n}= C_2^{1/\varepsilon} a_n^{-\frac 1{2(1+\varepsilon)}} n^{-\frac{1}{2+\varepsilon}}$,  the latter quantity tends to zero, because $a_n \ge 1/\log(1+n)$, and so we can apply  Inequality~1 and  Proposition~1 (relation~(12)) given  on pages~453 and~455 respectively in  Shorack and Wellner~\cite{shorack}. Then we obtain
\begin{equation}
\label{proof_17}
\begin{split}
\textbf{P}_2&\leq 2\exp\Bigl( -\frac {\lambda^2}{2p_k}\,\frac 1{1+2\lambda/(3p_k \sqrt{n})}\Bigr)\\
& = 2 \exp\Bigl( -\frac {1}{2p_k} C_2^{2/\varepsilon} a_n^{-\frac 1{1+\varepsilon}} n^{\frac {\varepsilon}{2+\varepsilon}} [1+o(1)]\Bigr).
\end{split}
\end{equation}
From~\eqref{proof_phi} and~\eqref{proof_17} it follows that
\begin{equation}
\label{proof_18}
\textbf{P}_2 = [1-\Phi(x)]o(1), \ \text{ as} \ n\to\infty,
\end{equation}
uniformly in the range $-A\le x\le a_n z_n$. So, the first relation in~\eqref{proof_12} follows directly from~\eqref{proof_13}, \eqref{proof_15} and~\eqref{proof_18}.

The next step we prove the second relation in~\eqref{proof_12}. We have
\begin{equation}
\label{proof_19}
\begin{split}
&\textbf{P}\{\sqrt{n}|\an-\alp| \emph{D}_{n,u}^{\varepsilon} > C\,\delta \}\\
\leq & \sum_{i:\,i/n\in \Gamma_n}  \textbf{P}\{\sqrt{n}|\an-\alp| |\uin-i/n|^{\varepsilon}\vee |\uin-(i-1)/n|^{\varepsilon} > C\,\delta \} .
\end{split}
\end{equation}
By condition {\bf (iii)}, there exists $M>0$ such that $|\an-\alp|\le M n^{-1/(2+\varepsilon)}$ for all sufficiently large $n$, hence each item of the sum on the l.h.s. in~\eqref{proof_19} does not exceed
\begin{equation}
\label{proof_19a}
\textbf{P}\{\sqrt{n} |U_{i:n}-i/n| >\lambda \} + \textbf{P}\{\sqrt{n} |U_{i:n}-(i-1)/n| > \lambda \}, \quad i/n\in \Gamma_n.
\end{equation}
where $\lambda=C_{\varepsilon} a_n^{-1/(2\varepsilon)} n^{\frac {\varepsilon}{2(2+\varepsilon)}}$ and     $C_{\varepsilon}=(C/M)^{1/\varepsilon}$. Obviously (cf.~\eqref{proof_16}-\eqref{proof_17}), it suffices to prove the desired bound for the first of two probabilities in~\eqref{proof_19a}. Applying once more the exponential Inequality~1 for uniform order statistics (cf. Shorack and Wellner~\cite{shorack}, pp.~453,~455) and the fact that $|i/n-\alp|\leq Mn^{-1/(2+\varepsilon)}$  for all sufficiently large $n$, we obtain
\begin{equation*}
\textbf{P}\{\sqrt{n} |U_{i:n}-i/n| >\lambda \} \leq 2\exp\lr -\frac {1}{2\alp} C_{\varepsilon}^2 a_n^{-\frac 1{\varepsilon}} n^{\frac {\varepsilon}{2+\varepsilon}} \Bigl[1+O(n^{-1/(2+\varepsilon)})\Bigr]\rr.
\end{equation*}
Since  the number of items on the r.h.s. in \eqref{proof_19} does not exceed  $|n|\alp-\an|+1|=O(n^{\frac{1+\varepsilon}{2+\varepsilon}})$, the latter bound implies  that the quantity on the r.h.s. in \eqref{proof_19} is of the order
\begin{equation*}
\label{proof_20_}
n^{\frac{1+\varepsilon}{2+\varepsilon}} \exp\Bigl( -\frac {1}{2\alp} C_{\varepsilon}^2 a_n^{-\frac 1{\varepsilon}} n^{\frac {\varepsilon}{2+\varepsilon}} \Bigl[1+o(1)\Bigr]\Bigr).
\end{equation*}
This together with~\eqref{proof_phi} imply the required relation.

It remains to prove the last relation in~\eqref{proof_12}. Fix some $\gamma>0$ such that $[\alp-\gamma,\alp+\gamma]\subseteq U_{\alp}$, set $r_n=k\wedge \kn$, $s_n=k\vee \kn+1$, where $\kn= n\an$ (cf.~\eqref{tn}). Then
\begin{equation}
\label{proof_20}
\textbf{P}\Bigl( \bigcup_{i:\,i/n\in \Gamma_n}\lf \uin \notin U_{\alp}\rf\Bigr)\leq \textbf{P}(U_{r_n:n} < \alp -\gamma) +
\textbf{P}(U_{U_{s_n:n}:n} > \alp+\gamma).
\end{equation}
Observe that both sequences $r_n/n$ and $s_n/n$ satisfy condition~{\bf (iii)}, along with the sequence $\an=\kn/n$.
Let us estimate the first probability on the r.h.s. in~\eqref{proof_20} (the treatment of
the second one is similar).
Define a~binomial random variable $S_n=\sharp \{i : U_i < \alp-\gamma \}$, then the first term on the r.h.s. in~\eqref{proof_20} is equal to
\begin{equation}
\label{proof_12_2}
\begin{split}
\textbf{P}(S_n \geq r_n) = &\textbf{P}\bigl(S_n-\textbf{E}S_n \geq r_n -n\alp  +\gamma n\bigr)\\
 = &\textbf{P}\bigl(n^{-1}(S_n-\textbf{E}S_n)\ge \gamma+ o(1)\bigr)
\end{split}
\end{equation}
and by the~classical Hoeffding~\cite{ho} inequality, the latter quantity is no greater than
$\exp(-2n(\gamma+ o(1))^2)$, which is $[1-\Phi(x)]o(1)$, uniformly in the range $-A \le x \le a_n n^{1/2}$,
and the last relation in~\eqref{proof_12} follows.

Relations~\eqref{proof_11}-\eqref{proof_11_} and~\eqref{proof_12} directly imply~\eqref{proof_10}, which yields~\eqref{proof_5}.

\noindent{\bf Proof of \eqref{proof_6}}. \ By condition~{\bf (iv)}, there exists $b>0$ such that
\begin{equation*}
\label{proof_21}
\sqrt{n}|V_n|\leq bn^{-\varepsilon/(2(2+\varepsilon))}(|X_{(\kn+1):n}|\vee |X_{(n-\mn):n}|),
\end{equation*}
for all sufficiently large $n$. Thus,
\begin{equation*}
\label{proof_22}
\textbf{P}\lr \sqrt{n}|V_n|/\sigma >\delta \rr \leq \textbf{P}\lr |X_{(\kn+1):n}|\vee |X_{(n-\mn):n}|>\sigma a_n^{-1/2}\rr\leq \textbf{P}_{3}+\textbf{P}_{4},
\end{equation*}
where $\textbf{P}_{3}=\textbf{P}\bigl(|X_{(\kn+1):n}|>\sigma a_n^{-1/2}\bigr) $,  \
$\textbf{P}_{4}=\textbf{P}\bigl(|X_{(n-\mn):n}|>\sigma a_n^{-1/2}\bigr) $. Let us estimate $\textbf{P}_{3}$ (the treatment for $\textbf{P}_{4}$ is same and therefore omitted). We have
\begin{equation}
\label{proof_23}
\begin{split}
\textbf{P}_{3}&=\textbf{P}\lr \lv F^{-1}(U_{(\kn+1):n})\rv >\sigma a_n^{-1/2}\rr \\
&\leq \textbf{P}\lr \lv F^{-1}(U_{(\kn+1):n})- F^{-1}(\alp)\rv +\lv F^{-1}(\alp)\rv  >\sigma a_n^{-1/2}\rr,\\
&\leq \textbf{P}\lr \lv U_{(\kn+1):n}- \alp\rv^{\varepsilon} >\sigma a_n^{-1/2} (1+o(1))\rr + \textbf{P}(U_{(\kn+1):n} \notin U_{\alp}).
\end{split}
\end{equation}
Observe that the first term on the r.h.s. in~\eqref{proof_23} is equal to zero, for all sufficiently large $n$, and the second one
is $[1-\Phi(x)]o(1)$, uniformly in the range $-A\le x\le a_n z_n$.  This completes the proof of~\eqref{proof_6} and the theorem. \qed

\medskip
\noindent{\bf Proof of Theorem~\ref{thm2}}. Let us first prove relation~\eqref{thm_2}. By Lemma~\ref{lem_2.1} and relation~\eqref{proof_1}, we have
\begin{equation*}
\label{proof t21}
\text{\em Var}(\tn)=\text{\em Var}(\widetilde{L}_n)+\text{\em Var}(R_n+V_n)+2 \text{\em cov}(\widetilde{L}_n,R_n+V_n).
\end{equation*}
Since $W_i$ are bounded, all conditions of Theorem~2\,(ii)~\cite{vv82} are satisfied, and hence
\begin{equation*}
\label{proof t22}
\sigma^{-1} n^{1/2}\sqrt{\text{\em Var}(\widetilde{L}_n)}=1+O(n^{-1/2})
\end{equation*}
(cf. \cite{vv82}, p.~431). \ Furthermore, we have
\begin{equation*}
\label{proof t21}
\begin{split}
n|\text{\em cov}(\widetilde{L}_n,R_n+V_n)|&\leq n[\text{\em Var}(\widetilde{L}_n)\text{\em Var}(R_n+V_n)]^{1/2}\\
&= \sigma[n\text{\em Var}(R_n+V_n)]^{1/2}(1+O(n^{-1/2})).
\end{split}
\end{equation*}
The latter three relations imply that in order to prove~\eqref{thm_2}, it suffices to show that
\begin{equation}
\label{proof t23}
n \text{\em Var}(R_n+V_n)=O\bigl(n^{-\frac{2\varepsilon}{2+\varepsilon}}\bigr).
\end{equation}
We have
\begin{equation}
\label{proof t24}
n \text{\em Var}(R_n+V_n)\leq n \textbf{E}(R_n+V_n)^2\leq  5n\Bigl[\sum_{i,j=1}^{2} \textbf{E}\bigl(I_j^{(i)}\bigr)^2 + \textbf{E}V_n^2 \,\Bigr],
\end{equation}
where $I_j^{(i)}$ are as in \eqref{proof_9}-\eqref{proof_10}. We will show that
\begin{equation}
\label{proof t25}
n \textbf{E}\bigl(I_j^{(1)}\bigr)^2 =O(n^{-\varepsilon})=o\bigl(n^{-\frac{2\varepsilon}{2+\varepsilon}}\bigr), \ \  n \textbf{E}\bigl(I_j^{(2)}\bigr)^2 =O(n^{-\frac{2\varepsilon}{2+\varepsilon}}),\ j=1,2,
\end{equation}
and that
\begin{equation}
\label{proof t26}
n \textbf{E}V_n^2 =O(n^{-\frac{2\varepsilon}{2+\varepsilon}}).
\end{equation}
Relations \eqref{proof t24}-\eqref{proof t26} imply the desired bound~\eqref{proof t23}.

We first prove~\eqref{proof t25}, and consider in detail only the case $j=1$ (the treatment in the case $j=2$ is same and therefore omitted). Let as before $k=[\alpha n]$ and $k_n=\alpha_n n$. By \eqref{proof_11} and the Schwarz inequality, we have
\begin{equation*}
\label{proof t27}
\begin{split}
\textbf{E}\bigl(I_1^{(1)}\bigr)^2 \leq & K^2 [\textbf{E}(A_n-\alp)^4\textbf{E}(X_{k:n}-\xia)^4]^{1/2}\\
= & K^2 n^{-2}[\textbf{E}(\na-\alp n)^4\textbf{E}(X_{k:n}-\xia)^4]^{1/2}.
\end{split}
\end{equation*}
By well-known formula for 4-th moments of a~binomial random variable, we have $\textbf{E}(\na-\alp n)^4=3\alp^2(1-\alp^2)n^2(1+o(1))$. Thus, there exists a~positive constant $C$ independent of $n$ such that
\begin{equation}
\label{proof t28}
n\textbf{E}\bigl(I_1^{(1)}\bigr)^2 \leq C [\textbf{E}(X_{k:n}-\xia)^4]^{1/2}
\end{equation}
for all sufficiently large $n$. We have
\begin{equation}
\label{proof t29}
\begin{split}
&\textbf{E}(X_{k:n}-\xia)^4= \textbf{E}(F^{-1}(U_{k:n})-F^{-1}(\alp))^4\\
=&\textbf{E}[(F^{-1}(U_{k:n})-F^{-1}(\alp))^4 \textbf{1}_{\{U_{k:n}\in U_{\alp}\}}]\\
+& \textbf{E}[(F^{-1}(U_{k:n})-F^{-1}(\alp))^4 \textbf{1}_{\{U_{k:n}\notin U_{\alp}\}}]\\
\leq & C_H^4 \textbf{E}|U_{k:n}-\alp|^{4\varepsilon}+ \textbf{E}[(X_{k:n}-\xia)^6]^{2/3}[\textbf{P}(U_{k:n}\notin U_{\alp})]^{1/3},
\end{split}
\end{equation}
where $C_H$ is a~constant from the H\"{o}lder condition {\bf (ii)}. Note that if $\varepsilon>1/2$, then
$\textbf{E}|U_{k:n}-\alp|^{4\varepsilon}\leq \textbf{E}|U_{k:n}-\alp|^{2}=O(n^{-1})$, and if $\varepsilon \leq 1/2$,
then $\textbf{E}|U_{k:n}-\alp|^{4\varepsilon}\leq (\textbf{E}|U_{k:n}-\alp|^{2})^{2\varepsilon}=O(n^{-2 \varepsilon})$.
Since  moments of any order of $X_{k:n}$ are finite for all sufficiently large $n$ and because $\textbf{P}(U_{k:n}\notin U_{\alp})=O(\exp(-cn))$ with some $c>0$ (cf.~\eqref{proof_20}-\eqref{proof_12_2}), the latter bounds and relations~\eqref{proof t28}-\eqref{proof t29} imply the first of relations~\eqref{proof t25}.

Consider $I_1^{(2)}$. By condition {\bf (iii')}, there exists $d>0$ such that \\ $(\an-\alp)^2\leq d n^{-(2+\varepsilon-\varepsilon^2)/(2+\varepsilon)}(\log n)^{-\varepsilon}$, for all sufficiently large $n$. Then in view of~\eqref{proof_11_} we obtain
\begin{equation}
\label{proof t210}
n\textbf{E}\bigl(I_1^{(2)}\bigr)^2 \leq n K^2 (\aln-\alp)^2 \textbf{E} \emph{D}_n^2 \leq K^2n^{\frac{\varepsilon^2}{2+\varepsilon}}(\log n)^{-\varepsilon}\textbf{E} \emph{D}_n^2.
\end{equation}
Hence, to get the second bound in~\eqref{proof t25}, it suffices to show that
\begin{equation}
\label{proof t211}
\textbf{E}\emph{D}_n^2=O\lr (\log n)^{\varepsilon}n^{-\varepsilon}\rr.
\end{equation}
For all sufficiently large $n$, \ $\an \in U_{\alp}$ and
\begin{equation*}
\label{proof t212}
\textbf{E}\emph{D}_n^2\leq C_H^2 \textbf{E}\lr\emph{D}_{n,u}^{\varepsilon}\rr^2=
C_H^2\textbf{E}\lr \max_{i:\,i/n\in \Gamma_n}|\uin-i/n|^{2\varepsilon}\vee |\uin-(i-1)/n|^{2\varepsilon}\rr,
\end{equation*}
where $\emph{D}_n$ is as in~\eqref{proof_12a}. The latter quantity does not exceed
\begin{equation}
\label{proof t213}
\begin{split}
&t^{\varepsilon}C_H^2 (\log n)^{\varepsilon}n^{-\varepsilon} \\
+ &\textbf{P}
\lr \bigcup_{i:\,i/n\in \Gamma_n} \lf |\uin-i/n|\vee |\uin-(i-1)/n|> \sqrt{t\frac{\log n}n}\rf\rr \\
\leq & t^{\varepsilon}C_H^2 (\log n)^{\varepsilon}n^{-\varepsilon}+|\alp n-\kn+1|\lr \textbf{P}_1+\textbf{P}_2\rr,
\end{split}
\end{equation}
where $t$ is a constant which will be chosen later, and
\begin{equation*}
\label{proof t214}
\textbf{P}_1= \textbf{P}\lr |\uin-i/n|> \sqrt{t\frac{\log n}n}\rr,\ \ \textbf{P}_2=\textbf{P}\lr |\uin-(i-1)/n|> \sqrt{t\frac{\log n}n}\rr.
\end{equation*}
It is obvious that both $\textbf{P}_1$ and $\textbf{P}_2$ are of the same order of magnitude,  so it suffices to estimate $\textbf{P}_1$, where we can apply once more the  Inequality~1 from Shorack and Wellner~\cite{shorack}. We have
\begin{equation*}
\label{proof t215}
\textbf{P}_1= \textbf{P}\lr \sqrt{n}|\uin-i/n|> \sqrt{t\log n}\rr\leq 2\exp\lr-\frac{t\log n}{2\alp}(1+O(|\an-\alp|))\rr,
\end{equation*}
hence if we choose $t\ge 4\alp$, we obtain $\textbf{P}_1+\textbf{P}_2=O(n^{-2})$, and the second term on
the r.h.s. in~\eqref{proof t213} becomes negligible in order relative to the first one.
This proves~\eqref{proof t211} and the second relation in~\eqref{proof t25}.

We now turn to the proof of~\eqref{proof t26}. By condition {\bf (iv')}, the exists a~constant $C>0$, not depending on
$n$, such that  $\sum_{i=k_n+1}^{n-m_n}|\cin -\cin^0|\leq Cn^{\frac{2-\varepsilon}{2(2+\varepsilon)}}$, for all sufficiently large $n$, and
\begin{equation*}
\begin{split}
n\textbf{E}V_n^2 &\leq n^{-1} \lr \sum_{i=k_n+1}^{n-m_n}|\cin -\cin^0|\rr^2\textbf{E}\bigl(X^2_{k_n+1:n}\vee X^2_{n-m_n:n}\bigr)\\
&\leq C^2n^{-1} n^{\frac{2-\varepsilon}{2+\varepsilon}} \textbf{E}\bigl(X^2_{k_n+1:n}\vee X^2_{n-m_n:n}\bigr)=O\bigl( n^{-\frac{2\varepsilon}{2+\varepsilon}}\bigr),
\end{split}
\end{equation*}
and~\eqref{proof t26} follows.

Thus, relation~\eqref{thm_2} is proved, and we are now in a~position to prove that relations~\eqref{thm_1} hold true if we replace $\sigma/n^{1/2}$  by $\sqrt{{\text{\em Var}}(\tn)}$. We prove the first of relations~\eqref{thm_1},  the second one will then follow from the first if we replace $\cin$ by $-\cin$.

Fix an~arbitrary sequence $a_n\to 0$ and $A>0$, set $\lambda_n=\sigma^{-1}n^{1/2} \sqrt{{\text{\em Var}}(\tn)}$ and write
\begin{equation}
\label{proof t217}
\frac{\textbf{P}\bigl((\tn-\mu_n)/ \sqrt{{\text{\em Var}}(\tn)} >x\bigr)}{1-\Phi(x) }=\frac {1-F_{\tn}(\lambda_nx)}{1-\Phi(\lambda_n x)} \, \,\frac{1-\Phi(\lambda_n x)}{1-\Phi(x)}.
\end{equation}
Set $B=A\sup_{n\in \mathbb{N}}\lambda_n$ and $b_n=\lambda_n a_n$, Since $\lambda_n \to 1$, the  number $B$ exists and $b_n \to 0$. Hence, by Theorem~\ref{thm1}, the first ratio on the r.h.s. in ~\eqref{proof t217} tends to $1$ as $\nty$, uniformly in $x$ such that $-B \leq \lambda_n x \leq b_n z_n$, where $z_n=n^{\varepsilon/(2(2+\varepsilon))}$, in particular,  uniformly in the range $-A \leq  x \leq a_n z_n$. Furthermore, we see that $|\lambda_n-1|^{1/2} a_nz_n \to 0$, which is due to the fact that $|\lambda_n-1|^{1/2}=O\bigl(n^{-\frac {\varepsilon}{2(2+\varepsilon)}}\bigr)$.  Hence, by Lemma A1 from Vandemaele and Veraverbeke~\cite{vv82},  the second ratio on the r.h.s. in ~\eqref{proof t217} also tends to~$1$, uniformly in the range $-A \leq  x \leq a_n z_n$. The theorem is proved.  \qed

{\bf Acknowledgments.} \ The author is grateful to the~referee for his valuable remarks and suggestions that led to  improvement of the article.

\end{document}